\documentclass[12pt]{amsart}
\usepackage{amssymb}

 \rm

\renewcommand{\(}{\left(}
\renewcommand{\)}{\right)}
\renewcommand{\[}{\left[}
\renewcommand{\]}{\right]}

\newcommand{\x}{\times}

\newcommand{\abs}[1]{\left\lvert#1\right\rvert}
\newcommand{\norm}[1]{\left\lVert#1\right\rVert}
\newcommand{\st}{\:|\:}

\renewcommand{\phi}{\varphi}

\renewcommand{\Re}{{\mathrm{Re}\,}}

\newcommand{\CC}{\mathbb{C}}
\newcommand{\RR}{\mathbb{R}}

\newcommand{\I}{\mathrm{I}}
\newcommand{\II}{\mathrm{II}}
\newcommand{\Res}{\mathrm{Res}}

\theoremstyle{plain}
\newtheorem{thm}{Theorem}[section]
\newtheorem{lem}[thm]{Lemma}

\newtheorem{ex}[thm]{Example}

\theoremstyle{definition}
\newtheorem{defn}[thm]{Definition}

\theoremstyle{remark}

\title{The Brylinski beta function of a surface}

\author{E.~J.~ Fuller}
\address{Department of Mathematics, West Virginia University,
Morgantown, WV 26501}

\author{M.~K.~Vemuri}
\address{Chennai Mathematical Institute, Plot H1, SIPCOT IT Park,
Padur~PO, Siruseri 603 103.}

\dedicatory{For Baxter and Tonu.}

\date{September 25, 2010}

\begin{document}

\begin{abstract}
An analogue of Brylinski's knot beta function is defined for a submanifold
of $d$-dimensional Euclidean space.  This is a meromorphic function on the
complex plane.  The first few residues are computed for a surface in three
dimensional space.
\end{abstract}

\maketitle


\section{Introduction}\label{S:intro}

In \cite{bry}, Brylinski introduced the beta function of a geometric knot in
$\RR^3$.  He was partly motivated by the desire to give a definition of
M\"obius energy (see \cite{fre}) independent of an arbitrary
``renormalization''.  However, he also gave some beautiful formulae for the
first few residues of his beta function.  They turn out to be integrals of
polynomials in the curvature, torsion and their derivatives.

In this note, we consider arbitrary submanifolds of $\RR^d$.  Essentially the
same definition (as in \cite{bry}) works in this situation, to define the
beta function on the right half-plane, and it can be analytically continued
to be a meromorphic function on $\CC$, with only simple poles.  The location
of the poles is dependent on the dimension of the submanifold, and if $M$
is a hypersurface, the
residues are integrals of polynomials in complete contractions of the covariant derivatives
of the second fundamental form.

We consider surfaces in $\RR^3$ in more detail, and compute some residues.
In particular, we characterize the spheres by the vanishing of the residue
at $s=-4$.

For a knot, the beta function does not have a pole at the M\"obius invariant
parameter $s=-2$, and its value there coincides with the M\"obius energy.
However, for a surface, the M\"obius invariant parameter is $s=-4$, and the
beta function has a pole there, in general.  The value obtained after
subtracting the pole may be considered the natural renormalization of
M\"obius energy (see \cite{kus}).

\section{The Brylinski beta function}

Let $M$ be a compact smooth $n$-dimensional submanifold of $\RR^d$.  Let $dA$ denote the $n$-dimensional
area element
of $M$.  Observe that if $\Re s > -n$, and $F(u,v) = \norm{v-u}^s$, then
$F \in L^1(M \x M)$.

\begin{defn}
The {\em Brylinski beta function of $M$} is the function $B_M(s)$, defined
for $\Re s > -n$ by
$$
B_M(s) = \int_{M \x M} \norm{v-u}^s \, dA(v) \, dA(u).
$$
\end{defn}

We also need to consider the pointwise version of the beta function.  For fixed $u \in M$,
we define
$$
B_M^u(s) = \int_{M} \norm{v-u}^s \, dA(v), \qquad \Re s > -n.
$$
Thus
$$
B_M(s) = \int_{M} B_M^u(s) \, dA(u).
$$

Note that $B_M^u$ and $B_M$ are analytic in the half-plane $\Re s > -n$.

\begin{ex}
Let $M=S^2(r):=\{ u \in \RR^3 \st \norm{u}=r \}$ be the sphere of radius $r$ in $\RR^3$.
Then $B_M^u(s) = \frac{2^{s+3} \pi r^{s+2}}{s+2}$ for all $u \in M$, and so
$$
B_{S^2(r)}(s) = \frac{2^{s+5} \pi^2 r^{s+4}}{s+2}
$$
\end{ex}

\begin{proof}
By rotational invariance, it is clear that $B_M^u(s)$ is independent of $u$, so we
take $u=(0,0,r)$ for the computation.  We use Cartesian coordinates $(x, y, z)$
and spherical coordinates $(r, \theta, \phi)$, thus
$$
\begin{aligned}
x &= r \sin \theta \cos \phi \\
y &= r \sin \theta \sin \phi \\
z &= r \cos \theta.
\end{aligned}
$$
So
$$
\begin{aligned}
\norm{v-u}^s
= & \; (x^2 + y^2 + (z - r)^2)^{s/2} \\
= & \; \[r^2 \sin^2 \theta \cos^2 \phi + r^2 \sin^2 \theta \sin^2 \phi +
            r^2 (\cos \theta - 1)^2\]^{s/2} \\
= & \; \[ r^2 \sin^2 \theta + r^2 (\cos^2 \theta - 2 \cos \theta + 1) \]^{s/2} \\
= & \; r^s (2 - 2 \cos \theta)^{s/2} \\
= & \; 2^s r^s \sin^s(\theta/2).            
\end{aligned}
$$

We recall that the area element on  $M$ is $dA = r^2 \sin \theta \, d\theta \, d\phi$, so
$$
\begin{aligned}
B_M^u(s)
= & \; \int_M \norm{v-u}^s \, dA(v) \\
= & \; 2^s r^s \int_0^{2\pi} \int_0^\pi \sin^s(\theta/2) \, r^2 \sin(\theta) \, d\theta \, d\phi \\
= & \; 2^{s+2} \pi r^{s+2} \int_0^\pi \sin^{s+1}(\theta/2) \cos(\theta/2) \, d\theta \\
= & \; 2^{s+3} \pi r^{s+2} \int_0^1 t^{s+1} \, dt \\
= & \; \frac{2^{s+3} \pi r^{s+2}}{s+2}.
\end{aligned}
$$
\end{proof}

\begin{ex}
More generally, if $M=S^n(r):=\{ u \in \RR^{n+1} \st \norm{u}=r \}$ is the sphere of radius
$r$ in $\RR^{n+1}$, then
$$
B_M^u(s) = 2^{s+n} \omega_{n-1} r^{s+n} B(\frac{s+n}{2}, \frac{n}{2})
$$ for all $u \in M$,
and so
$$
B_{S^n(r)}(s) = 2^{s+n} \omega_{n-1} \omega_n r^{s+2n} B(\frac{s+n}{2}, \frac{n}{2}),
$$
where $\omega_n$ denotes the $n$-dimensional ``area" of $S^n(1)$, and
$B(s, t)$ is Euler's beta function.
\end{ex}

\begin{proof}
$$
\begin{aligned}
B_M^u(s)
= & \; \int_0^\pi \[2r \sin(\theta/2)\]^s \, \omega_{n-1} (r \sin \theta)^{n-1} r \, d\theta \\
= & \; 2^s \omega_{n-1} r^{s+n} \int_0^\pi \sin^s(\theta/2) \sin^{n-1}(\theta) \, d\theta \\
= & \; 2^{s+n-1} \omega_{n-1} r^{s+n} \int_0^\pi \sin^{s+n-1}(\theta/2)
                                                                            \cos^{n-1}(\theta/2) \, d\theta \\
= & \; 2^{s+n} \omega_{n-1} r^{s+n} B(\frac{s+n}{2}, \frac{n}{2}).                                                                         
\end{aligned}
$$
\end{proof}

\section{The analytic continuation}
We begin with two analytic lemmas that are used in the main argument.
We try to imitate the arguments in \S 3.2 and \S 3.9 of \cite{GS}.  There are two obstacles to this.
Firstly, in our setting, the test-function is also varying (holomorphically); this is easily overcome
using Lemma \ref{L:entire}.  The second obstacle is somewhat more serious: the test function is
actually not smooth at the origin.  This problem is resolved by ``blowing-up" the origin, and applying
the Malgrange preparation theorem to show that the pulled back test function extends smoothly
across the exceptional divisor (Lemma \ref{L:prep})

\begin{lem}\label{L:entire}
Let $X$ be a compact Hausdorff space.  Let $\mu$ be a finite Baire measure on $X$.  Suppose
$G:X \x \CC \to \CC$ is continuous, and for each $x \in X$, the function $G(x, \cdot)$ is entire.  Then
$$
g(s) = \int_X G(x, s) \, d\mu(x)
$$
is entire.
\end{lem}

\begin{proof}
Fix $R \in (0, \infty)$ and put $C=\sup\{G(x, s) \st x \in X, \, \abs{s}=R\}$.  By the Cauchy estimates,
$$
\abs{\frac{\partial^k G}{\partial s^k}(x, 0)} \le \frac{k! C}{R^k},
$$ 
so
$$
\abs{\frac{d^kg}{d s^k}(0)} = \abs{\int_X \frac{\partial^k G}{\partial s^k}(x, 0) \, d\mu(x)}  \le \frac{k! C \mu(X)}{R^k},
$$
and so the radius of convergence of the Taylor series of $g$ is at least $R$.  Since $R$ is arbitrary,
the result follows.
\end{proof}

Suppose $f \in C^\infty(\RR^n)$ and $f$ vanishes to second order at $0$.  If $n>1$, the function
$F(w)=f(w)/\norm{w}^2$ does not necessarily extend smoothly to $0$.  However, the singularity is mild,
and may be resolved by ``blowing up".  Define $P: \RR \x S^{n-1}(1) \to \RR^n$ by
$P(r, w) = rw$.

\begin{lem}\label{L:prep}
The function $F \circ P$ extends to a smooth function.
\end{lem}

\begin{proof}
Let $\Phi:\RR^{n-1} \to S^{n-1}(1)$ be a chart and define $g:\RR \x \RR^{n-1} \to \RR$ by
$$
g(r, w') = (f \circ P)(r, \Phi(w')) = f(r \Phi(w')).
$$
By the Malgrange preparation theorem (\cite[Theorem 7.5.6]{H}), we can write
$$
g(r, w') = q(r, w') r^2 + r_1(w') r + r_0(w'), \qquad q \in C^\infty(\RR^n), r_1, r_0 \in C^\infty(\RR^{n-1}).
$$
However, since $g$ vanishes to second order on the hyperplane $r=0$, the functions $r_1$ and $r_2$ are
identically zero, and we have $g(r,w')=q(r,w') r^2$.  But $(F \circ P)(r, \Phi(w'))$ agrees with $q(r,w')$ when
$r \ne 0$, so we can use $q$ to locally extend $F \circ P$.
\end{proof}

\begin{thm}\label{T:AC}
The function $B_M^u$ can be analytically continued to a meromorphic function on $\CC$
with simple poles at $-n-j$, $j=0,2,4,\dots$.  Moreover, if $M$ is a hypersurface, the residues are polynomials
in complete contractions of the covariant derivatives of the second fundamental form.
\end{thm}

\begin{proof}
If $\psi \in C_c^\infty(\RR^d)$ is identically $1$ in a neighborhood
of $u$ then the localized beta function
$$
B_M^\psi(s) = \int_{M} \norm{v-u}^s  \, \psi(v) \, dA(v)
$$
has the same principal part as $B_M^u$ because their difference extends to a holomorphic function on $\CC$.
So it suffices to prove the result with $B_M^u$ replaced by $B_M^\psi$ for an appropriate $\psi$.  By rotating
and translating $M$, we may assume that $u=0$ and the tangent space to $M$ is $\RR^n \subseteq \RR^d$
(clearly this process does not affect the beta functions).
Then, in a neighborhood of $0$,  $M$ is the graph of a function $f:\RR^n \to \RR^{d-n}$ which vanishes to
second order at $0$.  By making the neighborhood smaller, we may assume $\norm{f(w)} < \norm{w}$.
Choose $\psi$ to have support in this neighborhood.  Then
$$
B_M^\psi(s) = \int_{\RR^n} \(\norm{w}^2 + \norm{f(w)}^2\)^{s/2} \, \phi(w) \, dw
$$
where $\phi(w)=\psi(w, f(w)) A(w)$ and $A(w)$ is the area-density (it may be expressed in terms of the partial
derivatives of $f$).

Now,
$$
\begin{aligned}
B_M^\psi(s)
= & \; \int_{\RR^n} \norm{w}^s \(1 + \frac{\norm{f(w)}^2}{\norm{w}^2}\)^{s/2} \, \phi(w) \, dw\\
= & \; \int_0^\infty r^{s+n-1} S(r,s) \, dr
\end{aligned}
$$
where
$$
S(r,s) = \int_{S^{n-1}(1)} \(1 + \frac{\norm{f(rw)}^2}{\norm{rw}^2}\)^{s/2} \, \phi(w) \, d\sigma(w)
$$
and $\sigma$ is the surface measure on $S^{n-1}(1)$.  Note that $\frac{S(r,s)}{\omega_{n-1}}$
is the mean value of the function $\(1 + \frac{\norm{f(rw)}^2}{\norm{rw}^2}\)^{s/2} \, \phi(w)$ on the sphere of
radius $r>0$.  We can extend the definition of $S(r,s)$ to all real values of $r$ by writing
$$
S(r,s) = \int_{S^{n-1}(1)} G(r, w, s) \, d\sigma(w)
$$
where
$$
G(r,w,s) = \(1 + \frac{\norm{(f \circ P)(r,w)}^2}{r^2}\)^{s/2} \, (\phi \circ P)(r, w), \quad r \in \RR, \, s \in \CC. 
$$

By Lemma \ref{L:prep}, $G \in C^\infty(\RR \x S^{n-1}(0) \x \CC)$.  Moreover, for each $r \in \RR$ and
$w \in S^{n-1}(0)$, the function $G(r, w,\cdot)$ is entire, so by Lemma \ref{L:entire},
$S \in C^\infty(\RR \x \CC)$ and $S(r, \cdot)$ is entire for each $r \in \RR$.  Note that by the equality
of mixed-partials, the functions $\frac{\partial^j S}{\partial r^j}(r, \cdot)$ are also entire.  Since
$G(-r, -w, s) = G(r,w,s)$, it follows that $S(\cdot, s)$ is even, and so
$\frac{\partial^j S}{\partial r^j}(0, \cdot) = 0$ for odd $j$.

Now fix a positive integer $k$.  For $\Re s > -n$, we have
$$
\begin{aligned}
B_M^\psi(s)
= & \; \int_0^1 r^{s+n-1} \[ S(r, s) - S(0,s) - r \frac{\partial S}{\partial r}(0,s) - \cdots -
                       \frac{r^{k-1}}{(k-1)!} \frac{\partial^{k-1} S}{\partial r^{k-1}}(0,s) \] \, dr \\
   & \; + \int_1^\infty r^{s+n-1} S(r,s) \, dr +
            \sum_{j=0}^{k-1} \frac{1}{j! (s+n+j)} \frac{\partial^j S}{\partial r^j}(0,s).
\end{aligned}
$$
By Taylor's theorem, the first integral on the right is defined and holomorphic as a function of $s$ for
$\Re s > -n-k$, so the right hand side is a meromorphic function with only simple poles at $-n-j$ on the half-plane
$\Re s > -n-k$.  By our remark about the odd-order partial derivatives of $S$, it follows that $B_M^\psi$
does not actually have a pole at $-n-j$ for odd $j$.  Since $k$ is arbitrary,
this provides the desired analytic continuation of $B_M^\psi$ to $\CC$.
Observe that
$$
\Res_{s=-n-j} B_M^u = \frac{1}{j!} \frac{\partial^j S}{\partial r^j}(0, -n-j).
$$

Now suppose $M$ is a hypersurface, i.e. $d=n+1$.  Then the area density $A(w) = (1+\norm{\nabla f(w)}^2)^{1/2}$,
so for small $r$, the spherical mean
$$
S(r,s)=\int_{S^{n-1}(1)} \(1+\frac{f(rw)^2}{r^2}\)^{s/2} (1+\norm{\nabla f(w)}^2)^{1/2} d\sigma(w),
$$
and so $\frac{\partial^j S}{\partial r^j}(0,s)$ may be expressed as a polynomial in the Taylor coefficients of $f$
at $0$ and the moment integrals $\int_{S^{n-1}(1)} w^\alpha \, d\sigma(w)$.

Let $\mathrm{I}$ and $\mathrm{II}$ denote the first and second fundamental forms of $M$.  Using the local
parametrization $(w, f(w))$ of $M$, we find
$$
\mathrm{I}_{ij} = \delta_{ij} + \partial_i f \partial_j f, \quad\text{and}\quad
\mathrm{II}_{ij} = \frac{\partial_i \partial_j f}{(1+\norm{\nabla f}^2)^{1/2}}.
$$
Now, if $\alpha$ is a multi-index, the covariant derivative of the second fundamental form
$$
(\mathrm{II}_{ij;\alpha})(0) = (\partial^\alpha \partial_i \partial_j f)(0) + Q
$$
where $Q$ is a polynomial in the Taylor coefficients (at $0$) of $f$ of order less than or equal to
$\abs{\alpha}+1$.  Such a ``triangular" relation may be inverted to express the Taylor coefficients
of $f$ as polynomials in $(\mathrm{II}_{ij;\alpha})(0)$.

It follows that $\Res_{s=-n-j} B_M^u$ may be expressed as a polynomial in  $(\mathrm{II}_{ij;\alpha})(0)$.
However, $B_M^u$ is independent of the choice of $f$, so this polynomial is $\mathrm{O}(n-1)$-invariant.
By Weyl's {\em First Fundamental Theorem on Orthogonal Invariants}, it follows that the polynomial may
be re-expressed as a polynomial in the complete contractions of $\mathrm{II}_{ij;\alpha}(0)$.
(cf. \cite{Weyl}.  Also, \cite[\S4.2 and \S4.6]{Gray} give a modern exposition of the application
of Weyl's theorem in Geometry)
\end{proof}

\section{Surfaces in $\RR^3$}\label{S:surfaces}

In this section, assume $M \subseteq \RR^3$ is a surface.   We will compute
$\mathrm{Res}_{s=k} B_M(s)$ for $k=-2$, $-4$ and $-6$
in terms of $\I$, $\II$ and the first two covariant derivatives of $\mathrm{II}$.

Fix $u \in M$ and choose coordinates such that $\I_{ij} = \delta_{ij}$ at $x$.  Let $\II_{ij;k}$ and $\II_{ij;kl}$
denote the components of the first and second covariant derivatives of $\II$.  Define, at $u$
$$
\begin{aligned}
H_0 &= \II_{ii} \\
H_1 &= \II_{ij} \II_{ij} \\
H_2 &= \II_{ij;k} \II_{ij;k} \\
H_3 &= \II_{ii;jj} \\
H_4 &= \II_{ij} \II_{ij;kk} \\
H_5 &= \II_{ij} \II_{kk;ij}.
\end{aligned}
$$
Here we are using the extended Einstein summation convention, where we sum over repeated indices,
even if they are both covariant.  This is justified because we are working in a coordinate system which is
orthogonal at $u$.  One can get formulae for the $H_\alpha$ in a general coordinate system by first raising
one of the repeated indices using $\I$ before summing.

Being complete contractions of tensors, $H_\alpha$, $\alpha=0, \dots, 5$ are smooth functions on $M$,
independent of any coordinate choices.  Note that $H_0$ is just the mean-curvature.

\begin{thm}
$$
\begin{aligned}
\Res_{s=-2} B_M^u &= 2\pi \\
\Res_{s=-4} B_M^u &= \frac{\pi}{2} (2H_1 - H_0^2)\\
\Res_{s=-6} B_M^u &= \frac{\pi}{32} \(
-\frac{3}{16} H_0^4 + \frac{3}{4} H_1^2 + \frac{4}{3} H_2
-\frac{1}{2} H_0 H_3 + \frac{3}{2} H_4 + \frac{1}{2} H_5
\)
\end{aligned}
$$
\end{thm}

\begin{proof}
Assume, as in the proof of Theorem \ref{T:AC} that $u=0$ and the $xy$-plane is tangent to $M$ at $0$.  So
locally $M$ is the graph of a function $z=f(w)$, where $w=(x,y)$, which vanishes to second order at $0$.
By the proof of Theorem \ref{T:AC},
$$
\Res_{s=-2-j} B_M^u = \frac{1}{j!} \frac{\partial^j S}{\partial r^j}(0, -2-j).
$$
where
$$
S(r,s) = \int_{S^1(1)} 
\(1 + \frac{f(rw)^2}{r^2}\)^{s/2} \, (1 + \norm{\nabla f(rw)}^2)^{1/2} \, d\sigma(w),
$$
for small $r$

If we write
$$
\begin{gathered}
f(x, y) = \\
b_1 x^2 + b_2 x y + b_3 y^2 + \\
c_1 x^3 + c_2 x^2 y + c_3 x y^2 + c_4 y^3 + \\
d_1 x^4 + d_2 x^3 y + d_3 x^2 y^2 + d_4 x y^3 + d_5 y^4
+ O(\norm{w}^5),
\end{gathered}
$$
$$
\begin{gathered}
\(1 + \frac{f(rw)^2}{r^2}\)^{s/2} \, (1 + \norm{\nabla f(rw)}^2)^{1/2} \\
= 1 + \frac{s}{2} \frac{f^2}{r^2} + \frac12 \norm{\nabla f}^2 + \frac{s}{4} \frac{f^2 \norm{\nabla f}^2}{r^2}
         + \frac{s(s-2)}{8} \frac{f^4}{r^4} - \frac18 \norm{\nabla f}^4 + O(r^6),
\end{gathered}
$$
substitute into the definition of $S(r,s)$ and perform the indicated differentiation, we find
$$
\begin{aligned}
\Res_{s=-2} B_M^u &= 2\pi \\
\Res_{s=-4} B_M^u &= \frac{\pi}{2} (b_3^2 - 2 b_1 b_3 + b_2^2 + b_1^2) \\
\Res_{s=-6} B_M^u &= \frac{\pi}{32}
\(
\begin{gathered}
72 b_3 d_5 - 24 b_1 d_5 + 24 b_2 d_4 + 48 c_4^2 + 8 b_3 d_3 + 8 b_1 d_3 \\
 + 16 c_3^2 - 63 b_3^4 + 12 b_1 b_3^3 - 66 b_2^2 b_3^2 - 26 b_1^2 b_3^2 \\
- 44 b_1 b_2^2 b_3 - 24 d_1 b_3 + 12 b_1^3 b_3 + 24 b_2 d_2 + 16 c_2^2 \\
- 11 b_2^4 - 66 b_1^2 b_2^2 + 72 b_1 d_1 + 48 c_1^2 - 63 b_1^4
\end{gathered}
\)
\end{aligned}
$$
Using the coordinates $x, y$ on $M$, we can calculate the components of $\I$ and $\II$ to 
order two near $0$ and so determine the first and second covariant derivatives of $\II$ at $0$.  Using
this, we find that, at $0$, we have
$$
\begin{aligned}
H_0 =& \; 2 b_3+2 b_1 \\
H_1 =& \; 4 b_3^2+2 b_2^2+4 b_1^2 \\
H_2 =& \; 36 c_4^2+12 c_3^2+12 c_2^2+36 c_1^2 \\
H_3 =& \; 24 d_5+8 d_3-24 b_3^3-8 b_1 b_3^2-16 b_2^2 b_3-8 b_1^2 b_3-16 b_1 b_2^2+24 d_1-24 b_1^3 \\
H_4 =& \; 48 b_3 d_5+12 b_2 d_4+8 b_3 d_3+8 b_1 d_3-48 b_3^4-48 b_2^2 b_3^2-32 b_1^2 b_3^2 \\
     & -32 b_1 b_2^2 b_3+12 b_2 d_2-8 b_2^4-48 b_1^2 b_2^2+48 b_1 d_1-48 b_1^4 \\
H_5 =& \; 48 b_3 d_5+12 b_2 d_4+8 b_3 d_3+8 b_1 d_3-48 b_3^4-16 b_1 b_3^3-44 b_2^2 b_3^2 \\
     & -56 b_1 b_2^2 b_3-16 b_1^3 b_3+12 b_2 d_2-4 b_2^4-44 b_1^2 b_2^2+48 b_1 d_1-48 b_1^4,
\end{aligned}
$$
so the result follows.
\end{proof}

The quantity $2 H_1 - H_0^2$ equals $(p_1 - p_2)^2$ where $p_1$ and $p_2$ are the principal curvatures.
In particular, it is always non-negative and vanishes only at an umbilic point.  From this we conclude that
$\Res_{s=-4} B_M$ vanishes only for spheres.

It is well known that the ``warping" $(p_1 - p_2)^2 dA$ is
M\"obius invariant.  Of course, it depends only locally on $M$.  So the quantity
$$
\lim_{s \to -4} \(B_M(s) - \frac{\pi}{2(s+4)} \int_M (2 H_1 - H_0^2) \, dA\)
$$
may be thought of as the M\"obius energy of $M$ (see \cite{kus}).

\bibliographystyle{amsplain}
\bibliography{v9-bbs}

\end{document}